\documentclass[12pt]{article}
\usepackage{graphicx} 

\usepackage{cmap}  
\usepackage[T2A]{fontenc}
\usepackage[utf8]{inputenc}
\usepackage[english]{babel}

\usepackage[usenames]{color}
\usepackage{amsmath,amssymb,amsthm}
\usepackage[colorlinks=true,linkcolor=blue,urlcolor=red,unicode=true,hyperfootnotes=false,bookmarksnumbered]{hyperref}
\usepackage{mathtools}
\usepackage{xcolor}
\usepackage{csquotes}
\usepackage{wrapfig}
\usepackage{setspace}
\usepackage{enumitem}
\usepackage{appendix}
\usepackage{algorithm}
\usepackage{algpseudocode}
\usepackage{dsfont}
\usepackage{bm}
\usepackage{comment}

\newcommand{\ff}{\mathcal{F}}

\newcommand{\cG}{\mathcal{G}}

\renewcommand{\leq}{\leqslant}
\renewcommand{\geq}{\geqslant}
\renewcommand{\le}{\leqslant}
\renewcommand{\ge}{\geqslant}

\newtheorem{theorem}{Theorem}

\newtheorem{example}{Example}

\newtheorem{remark}{Remark}

\title{Short proofs of three combinatorial results in the Johnson scheme}

\author{Danila Cherkashin\thanks{Institute of Mathematics and Informatics, Bulgarian Academy of Sciences, Sofia, Bulgaria; E-mail: \url{jiocb@math.bas.bg}}, Yakov Shubin\thanks{Moscow Institute of Physics and Technology;
E-mail: \url{shubin.yakoff@gmail.com}}}

\begin{document}

\maketitle

\begin{abstract}
    In this note, we give short proofs of three theorems concerning extremal problems in the Johnson scheme, or, in other terminology, on \((n,k,L)\)-systems.

    The main result is a proof of the Aljohani--Bamberg--Cameron conjecture which claims that if $n > n_0(k)$ and there are an $(n,k,L)$-system and an $(n,k,\{0,\dots,k-1\}\setminus L)$-system whose sizes have product $\binom{n}{k}$, then they are a $t$-intersecting family and a Steiner system $S(t,k,n)$ for some $t$.
\end{abstract}

\section{Introduction}

For non-negative integers $n, a, b$, put $[n]=\{1, \ldots, n\}$ and, more generally, $[a, b]=\{a, a+1, \ldots, b\}$. Given a set $X$ and an integer $k \geqslant 0$, denote by $\binom{X}{k}$ the collection of all $k$-element subsets ($k$-sets) of $X$. A {\it family} is simply a collection of sets. 
For $L \subset [0,k-1]$ we say that $\ff \subset \binom{[n]}{k}$ is an \textit{$(n,k,L)$-system} if $|F_1 \cap F_2| \in L$ for any distinct $F_1,F_2 \in \ff$. The problem of determining the maximum size of an $(n,k,L)$-system is a classical problem in extremal set theory, introduced and studied by Deza, Erd\H{o}s and Frankl.

A 2-star is the family of all $k$-sets containing a fixed 2-element set. A family \(\mathcal F\subset\binom{[n]}{k}\) is called \(t\)-intersecting if
\(|F_1\cap F_2|\ge t\) for all distinct \(F_1,F_2\in\mathcal F\).

Let $G = (V,E)$ be a simple graph, denote by $\omega(G)$ its clique number, which is the size of a largest clique in $G$.
By $\alpha(G)$ we denote the size of a largest coclique (independent set), equivalently, the size of a largest clique in the complement of $G$.
A graph is called vertex-transitive if, given any two vertices $v_1$ and $v_2$ of $G$, there is an automorphism $f$ such that $f(v_1) = v_2$.
Denote by $e_U$ the number of edges induced by a subset $U \subset V(G)$.

Sometimes it is convenient to represent problems on $(n,k,L)$-systems in terms of the generalized Johnson graph $J(n,k,L)$. The vertex set of this graph is $\binom{[n]}{k}$, and two vertices are joined by an edge if the size of the intersection of the corresponding sets belongs to $L$. In these terms, the Deza--Erd\H{o}s--Frankl problem is to determine the clique number $\omega(J[n,k,L])$. 

These graphs (for specific $L$) are also studied as constant-weight codes. 
From the point of view of algebraic combinatorics, all graphs $J(n,k,L)$, for fixed $n$ and $k$, belong to the same Johnson association scheme and therefore have a common eigenbasis.

The graph $J(n,k,L)$ is vertex-transitive, and hence the clique-coclique bound
applies to it (see, for instance,~\cite{GM}):
\begin{equation}\label{cliq-co}
    \alpha(J[n,k,L]) \cdot \omega(J[n,k,L]) \le |V(J[n,k,L])| = \binom{n}{k}.
\end{equation}

In Section~\ref{sec2}, we show that, when $n$ is sufficiently large in terms of $k$, the clique-coclique bound~\eqref{cliq-co} can also be derived from the celebrated theorem of Deza, Erd\H{o}s and Frankl.

Aljohani, Bamberg and Cameron~\cite{ABC} studied when equality can occur in
\eqref{cliq-co}. In their terminology this means that the Johnson scheme $J(n, k)$ is non-separating.

There are several natural examples. 
\begin{example}[Aljohani--Bamberg--Cameron,~\cite{ABC}] \label{ex1}
Suppose that there exists a projective plane of order $q$, and put $n=q^2+q+1$ and $k=q+1$. The lines of the
projective plane form a clique in $J(n,k,\{1\})$, since any two lines intersect
in exactly one point. On the other hand, the family of all $k$-subsets containing
a fixed pair of points is a coclique in this graph. The sizes of these two families
multiply to $\binom{n}{k}$.
\end{example}


A more general source of examples comes from designs. Recall that a Steiner system
$S(t,k,n)$ is a family $\mathcal D\subset\binom{[n]}{k}$ such that every
$t$-element subset of $[n]$ is contained in exactly one member of $\mathcal D$. A projective plane of order $q$, used in Example~\ref{ex1}, is an example of a Steiner
system $S(2,q+1,q^2+q+1)$. In particular, any two distinct sets in $\mathcal D$ have intersection of size at most $t-1$. Thus, if $L=\{0,1,\ldots,t-1\}$, then the members of $\mathcal D = S(t,k,n)$ form a clique in $J(n,k,L)$. The family of all $k$-sets containing
a fixed $t$-element set is a coclique in $J(n,k,L)$, and we have
\[
    |\mathcal D|\binom{n-t}{k-t}
    =
    \frac{\binom{n}{t}}{\binom{k}{t}}\binom{n-t}{k-t}
    =
    \binom{n}{k}.
\]
The complementary case $L=\{t,t+1,\ldots,k-1\}$ gives the same construction with
the roles of the clique and the coclique interchanged. By the result of Keevash
on the existence of designs~\cite{keevash2014existence}, such Steiner systems
exist for all sufficiently large $n$ satisfying the necessary divisibility
conditions, namely
$$
\binom{k-i}{t-i}\mid \binom{n-i}{t-i}
\quad\text{for all } i=0,1,\ldots,t-1.
$$
See also~\cite{keevash2018existence, keevash2024short, glock2023existence}.

Aljohani, Bamberg and Cameron formulated the following conjecture. We prove it in
Section~\ref{sec2}. 

\begin{theorem} \label{mainth}
    Let $L$ be a non-empty proper subset of $\{0,\ldots,k-1\}$. If $\alpha(J[n,k,L]) \cdot \omega(J[n,k,L]) = \binom{n}{k}$ and $n$ is sufficiently large in terms of $k$, then there exists $t$ such that either $L$ or $\{0, \ldots, k-1\} \setminus L$ is equal to $\{0,1, \ldots, t-1\}$.
\end{theorem}

In~\cite{ABC}, this conjecture was proved for $k\le 4$. We also note that this conjecture appears as Problem~30.6 in the problem list from the 30th British Combinatorial Conference~\cite{BCC30}, 2024.

Now let us make a more detailed analysis of Example~\ref{ex1}. Projective planes of order $q$ exist for a prime power $q$; a major conjecture is that this does not hold for other values of $q$.
Example~\ref{ex1} and inequality~\eqref{cliq-co} give
\[
\alpha(J[q^2+q+1,q+1,\{1\}]) \leq \binom{q^2+q-1}{q-1}. 
\]
This inequality was proven for every $q$ by Cherkashin~\cite{cherkashin2024set} and then Linz~\cite{linz2026set} strengthened this result by proving that
\[
\alpha(J[n,k,\{1\}]) \leq \binom{n-2}{k-2} 
\]
for $3k-3 \leq n \leq k^2-k+1$.
Aljohani, Bamberg and Cameron also conjectured that every coclique of size $\binom{q^2+q-1}{q-1}$ in $J(q^2+q+1,q+1,\{1\})$ is a 2-star provided that $q > 2$; Cherkashin stated the same question for $J(k^2-k+1,k,\{1\})$ and arbitrary $k > 3$ (for $k=3$ there is a simple counterexample). 
Here we note that the Delsarte linear programming method (applied by Linz) together with the celebrated Ahlswede--Khachatrian Complete Intersection Theorem give the classification of all maximum cocliques in a graph $J(n,k,\{1\})$ for 
$3k-3 \leq n < k^2-k+1$.

\begin{theorem} \label{notmainth}
Let $I$ be a coclique in the graph $J(n,k,\{1\})$ of size $\binom{n-2}{k-2}$, $k \geq 3$. If $3k-3 < n < k^2 - k + 1$ then $I$ is a 2-star.
If $n = 3k-3$ then $I$ is either a 2-star or the union of all $k$-sets having at least 3 elements among fixed 4 elements.
\end{theorem}

Note that the most interesting case $n = k^2 - k + 1$ remains open. Finally, let us note that for the case of a large enough $k$ these results are covered by very general results of
Kupavskii and Zakharov~\cite{kupavskii2024spread} and Keller and Lifshitz~\cite{keller2021junta}; these results also cover the uniqueness issue.

After one determines the size of a maximum coclique and classifies all maximum cocliques in $G$, it is natural to consider stability and supersaturation problems. The results of Keller and Lifshitz were extended to the stability setting by 
Ellis, Keller and Lifshitz~\cite{ellis2024stability}.
The third result concerns the supersaturation problem, which is to determine the minimal number of edges in a vertex subset of a given size.
Many general tight results on this problem in the graphs $J(n,k,\{t\})$ were obtained in~\cite{kupavskii2026supersaturation}.

On the other hand, a complete asymptotic solution of this question is very complicated even in the case $J(n,3,\{1\})$, see~\cite{dubinin2024lower}. This problem has a threshold when the size of the vertex subset is \(cn^2\) for a constant $c$; in other regimes the supersaturation problem for the graph sequence $J(n,3,\{1\})$ is asymptotically solved.

\begin{theorem} \label{onemprenotmainth}
Let $c > 2/9$ be a constant and $U$ be a vertex subset of size $\lfloor cn^2\rfloor$ in $J(n,3,\{1\})$.
Then $U$ contains at least
\[
\left(\frac{9c^2}{2} - c + o(1) \right)n^3
\]
edges.
\end{theorem}

If $2c$ is an integer, then Theorem 9 from~\cite{dubinin2024lower} shows that Theorem~\ref{onemprenotmainth} is asymptotically tight. 
Also, the proof of Theorem~\ref{onemprenotmainth} gives tight lower bound in the case $n^2 = o(|U|)$, solved in \cite{shubin2020}.

Theorem~\ref{onemprenotmainth} can be rephrased as follows. Let $H$ be a 3-uniform hypergraph with $n$ vertices and $cn^2$ edges for some constant $c > 2/9$. Then $H$ has at least $\left(\frac{9c^2}{2} - c + o(1) \right)n^3$ pairs of edges with unit intersection.

\section{Proof of Theorem~\ref{mainth}}\label{sec2}

First, we state the condition of the celebrated theorem of Deza, Erd\H{o}s and Frankl, which gives an upper bound on the size of an $(n,k,L)$-system for $n > n_0(k)$. We then show how Theorem~\ref{mainth} follows from it. 

We will only need two parts of the theorem, so for convenience we state only those. For the proof and further information on $(n,k,L)$-systems, see~\cite{kupavskii2025delta}.

\begin{theorem}[Deza--Erd\H{o}s--Frankl,~\cite{DEF}]
\label{def}
    Fix a positive integer $k>0$ and a set 
    \[
    L=\left\{\ell_1<\ldots<\ell_r\right\} \subset[0, k-1].
    \]
    Assume that $n \geqslant n_0(k)$, and that $\mathcal{F}$ is an $(n, k, L)$-system.\\
(i) If $|\ff| \geqslant k^2 2^{r-1} n^{r-1}$ then
\[
(\ell_2-\ell_1) \mid (\ell_3-\ell_2) \mid \ldots \mid (\ell_r-\ell_{r-1}) \mid (k-\ell_r) .
\]
(ii) We have
\[
|\mathcal{F}| \leqslant \prod_{i=1}^r \frac{n-\ell_i}{k-\ell_i} .
\]
\end{theorem}

We now turn to our problem. Write
$L=\{\ell_1<\ldots<\ell_r\}\subset[0,k-1]$ and
$S=[0,k-1]\setminus L=\{s_1<\ldots<s_{k-r}\}$.
Let $\ff\subset\binom{[n]}{k}$ and $\cG\subset\binom{[n]}{k}$ be families
corresponding to a maximum clique and a maximum coclique in the graph
$J(n,k,L)$, respectively. Clearly, $\ff$ is an $(n,k,L)$-system, while $\cG$
is an $(n,k,S)$-system. By the assumption of Theorem~\ref{mainth}, we have
$|\ff|\cdot|\cG|=\binom{n}{k}$.

Assume that $n$ is larger than the constant $n_0$ from Theorem~\ref{def}.
By Theorem~\ref{def}, we have
\[
|\ff|\leqslant \prod_{i=1}^r \frac{n-\ell_i}{k-\ell_i} \qquad \mbox{and} \qquad 
|\cG|\leqslant \prod_{i=1}^{k-r} \frac{n-s_i}{k-s_i}.
\]
For the equality $|\ff|\cdot|\cG|=\binom{n}{k}$ to hold, both upper bounds must
be tight, since
\[
\prod_{i=1}^r \frac{n-\ell_i}{k-\ell_i}\cdot
\prod_{i=1}^{k-r} \frac{n-s_i}{k-s_i}
=
\binom{n}{k}.
\]

We may also assume that $|\ff|\geqslant k^2 2^{r-1} n^{r-1}$ and $|\cG|\geqslant k^2 2^{k-r-1} n^{k-r-1}$. Indeed, otherwise, using the upper bounds from part (ii) of Theorem~\ref{def},
we would have $|\ff|\cdot|\cG|\leqslant c(k)n^{k-1}$ for some constant $c(k)$, which is
impossible for sufficiently large $n$, since $|\ff|\cdot|\cG|=\binom{n}{k}$.

We know that $k-1$ belongs to either $L$ or $S$. First assume that $k-1\in L$.
Then $\ell_r=k-1$. By part (i) of Theorem~\ref{def}, we have
\[
(\ell_2-\ell_1)\mid(\ell_3-\ell_2)\mid\ldots\mid
(\ell_r-\ell_{r-1})\mid(k-\ell_r).
\]
Since $k-\ell_r=1$, all the differences
$\ell_2-\ell_1,\ldots,\ell_r-\ell_{r-1}$ are equal to $1$. Hence 
\[
L=\{k-r,k-r+1,\ldots,k-1\}.
\]

If, on the other hand, $k-1\in S$, then the same argument applied to the
$(n,k,S)$-system $\cG$ shows that $S=\{r,r+1,\ldots,k-1\}.$ Consequently, $L=\{0,1,\ldots,r-1\},$ again as required.

\section{Proof of Theorem~\ref{notmainth}}\label{sec3}

Define the Frankl family $\mathcal{F}_{n,k,t,r}$ by the following set of vertices in $J(n,k,\{t\})$
\[
\mathcal{F}_{n,k,t,r} := \left\{ S \in \binom{[n]}{k} : |S \cap [t+2r]| \geq t+r \right\}
\]
for some $n > k > t \geq 1$ and $r \geq 0$.

\begin{theorem}[Ahlswede--Khachatrian,~\cite{ahlswede1997complete}] \label{AKh}
A maximum $t$-intersecting family is always a Frankl set (up to a permutation of $[n]$).
Namely, $\mathcal{F}_{n,k,t,r}$ is the largest exactly in the range
\[
(k-t+1) \left (2+\frac{t-1}{r+1} \right) \leq n \leq (k-t+1) \left (2+\frac{t-1}{r} \right).
\]
(By convention $\frac{t-1}{r} = +\infty$ for $r=0$.)
\end{theorem}

Now let us return to the graph $J(n,k,\{1\})$.
Linz in~\cite{linz2026set} shows that the matrix $M = M(n,k)$ explicitly constructed by Wilson in~\cite{wilson1984exact} satisfies the following properties.
\begin{itemize}
\item[(i)] $M$ is a symmetric matrix of size $\binom{n}{k}$. There is a natural bijection $\pi$ between the rows/columns and the vertices of $J(n,k,\{1\})$ in which $M_{ij}$ depends only on $|\pi(i) \cap \pi(j)|$.\footnote{This means that $M$ belongs to the Bose--Mesner algebra of the Johnson scheme.}
\item[(ii)] For $3k-3 \leq n$ the largest eigenvalue of $M$ is $\binom{n-2}{k-2}$.
\item[(iii)] For $n \leq k^2-k+1$ the entries $M_{i,j}$ with $|\pi(i) \cap \pi(j)| \neq 1$ are at least 1.
\end{itemize}

The line-by-line repetition of the proof of item~(iii) gives that $M_{i,j}$ with $|\pi(i) \cap \pi(j)| = 0$ is strictly greater than 1 provided that $n$ is strictly smaller than $k^2 - k + 1$.
Now let $\chi_I$ be a characteristic vector of a coclique $I$ in the graph $J(n,k,\{1\})$ (according to the bijection $\pi$). Then in the range $3k-3 \leq n \leq k^2-k+1$ we have
\[
\binom{n-2}{k-2} |I| = \lambda_{max}(M)\cdot |I| \geq \langle M\chi_I,\chi_I \rangle \geq |I|^2,
\]
which is the Linz argument for $|I| \leq \binom{n-2}{k-2}$.

Now suppose we have $3k-3 \leq n < k^2-k+1$ and $|I| = \binom{n-2}{k-2}$. Then 
\[
\langle M\chi_I,\chi_I \rangle = |I|^2,
\]
and the mentioned specification of item~(iii) implies that a couple of sets from $I$ never have an empty intersection.
Thus $I$ is a 2-intersecting family and Theorem~\ref{AKh} finishes the proof.

\section{Proof of Theorem~\ref{onemprenotmainth}}\label{sec4}

Consider the adjacency matrix $A$ of the graph $J(n,3,\{1\})$. 
It has spectrum consisting of only four eigenvalues
\[
\lambda_0 = \frac{3(n-3)(n-4)}{2}, \quad \lambda_1 = \frac{(n-4)(n-9)}{2}, \quad \lambda_2 = -2n+11, \quad \lambda_3 = 3.
\]
(The calculation is standard and described for instance in~\cite{cherkashin2024set}.)
From now on, assume that $n \geq 7$, then the smallest eigenvalue is $\lambda_2$.

Let $U$ be a vertex subset of $J(n,3,\{1\})$, and let $\chi_{U}$ be the corresponding characteristic vector. 
Let $E_i$ be the eigenspace corresponding to $\lambda_i$, and consider the decomposition $\chi_U = \sum_{i=0}^3 \alpha_i u_i$, where $u_i \in E_i$ are unit vectors.
Then
\[
2e_U = \langle A \chi_U, \chi_U \rangle = \sum_{i=0}^3 \alpha_i^2 \lambda_i \geq \alpha_0^2 \lambda_0 + \sum_{i=1}^3 \alpha_i^2 \lambda_2.
\]
Note that $\sum_{i=0}^3 \alpha_i^2 = |U|$ which is the squared length of $\chi_U$. Also, since the graph is regular and connected, the eigenspace for $\lambda_0$ is one-dimensional and contains the all-ones vector. Thus $\alpha_0 = |U|/\sqrt{\binom{n}{3}}$. Putting everything together
\[
2e_U \geq (1+o(1)) \left ( \frac{3n^2}{2} \frac{|U|^2}{\binom{n}{3}}  - 2n \left(|U| - \frac{|U|^2}{\binom{n}{3}}\right) \right) = (1+o(1)) (9c^2 - 2c) n^3.
\]

\begin{remark}
The bound in Theorem~\ref{onemprenotmainth} can be obtained by using other matrices than just the adjacency matrix. 
On the other hand, it should be a linear programming (Delsarte--Schrijver) bound, and since it is asymptotically tight for infinitely many values of $c$, the result should be asymptotically the same.
\end{remark}

\paragraph{Acknowledgments.} Theorems~\ref{notmainth} and~\ref{onemprenotmainth}
are supported by Bulgarian NSF grant KP-06-N72/6-2023.

\bibliography{bibliography}
\bibliographystyle{plain}

\end{document}